\documentclass{amsart} \usepackage{amssymb,amsfonts,amscd}
\def\ni{{\noindent}}

\def\s01{{\{0,1\}}}
\def\cosec #1 {{\textrm {cosec}\ #1}}

\def\summ#1#2#3{{\sum_{#1=#2}^{#3}}}

\def\cupp#1#2#3{{\cup_{#1=#2}^{#3}}}

\def\prodd#1#2#3{{\prod_{#1=#2}^{#3}}}
\def\opluss#1#2#3{{\oplus_{#1=#2}^{#3}}}

\def\odott#1#2#3{{\odot_{#1=#2}^{#3}}}
\def\range #1#2#3{{{#1}={#2},\ldots,{#3}}}
\def\bx{{\bf x} }

\def\b1{{\bf 1} }

\def\del{{\delta}}
\def\Del{{\Delta}}

\def\alp{{\alpha}}
\def\veps{{\varepsilon}}
\def\bet{{\beta}}
\def\sig{{\sigma}}

\def\gam{{\gamma}}
\def\lam{{\lambda}}

\def\divides{{\mid}}
\def\notdivides{{\nmid}}
\def\bN{{\ensuremath{\mathbb N} }}
\def\bC{{\ensuremath{\mathbb C} }}
\def\bZ{{\ensuremath{\mathbb Z} }}

\def\cA{{\ensuremath{\cal A } }}
\def\cB{{\ensuremath{\cal B } }}
\def\cC{{\ensuremath{\cal C } }}

\def\cH{{\ensuremath{\cal H } }}

\def\cT{{\ensuremath{\cal T} }}

\def\bc{{\bf c}}
\def\bd{{\bf d}}
\def\be{{\bf e}}

\def\be{{\bf e}}
\def\b #1{{\bf #1}}

\def\proof{{\beginpf}}
\def\endproof{{\eopf}}

\def\lin{{\textrm{lin}}}
\def\linbar{{\overline{\textrm{lin}}}}

\newcommand\beq{\begin{equation}}
\newcommand\eeq{\end{equation}}
\newcommand\bdfn{\begin{defn}}
\newcommand\blem{\begin{lemma}}
\newcommand\elem{\end{lemma}}
\newcommand\bcor{\begin{cor}}
\newcommand\ecor{\end{cor}}
\newcommand\bthm{\begin{thm}}
\newcommand\ethm{\end{thm}}
\newcommand\edfn{\end{defn}}
\newcommand\bcas{\begin{cases}}
\newcommand\ecas{\end{cases}}

\def\nm#1{{\left\Vert #1 \right\Vert}}
\def\Nm{{\nm\cdot}}

\def\rad{\textrm {rad}}

\def\veps{{\varepsilon}}

\newcommand\sprod[2]{\langle{#1},{#2}\rangle}

\def\casif#1!{& \textrm{#1}\\}
\newif\ifrough
\roughfalse
\def\reff#1{{\ref {#1}}}
\def\l#1!{{\mbox{\label{#1}}}}

\def\refeq#1{{(\reff{#1})}}
\def\up#1{^{(#1)}}

\def\reflem#1{{Lemma \reff{#1}}}
\def\refdfn#1{{Definition \reff{#1}}}
\def\refthm#1{{Theorem \reff{#1}}}

\def\benum{\begin{enumerate}}
\def\eenum{\end{enumerate}}

\def\abws{{\bar \cA^{w*}}}

\newtheorem{thm}{Theorem}[section]
\newtheorem{defn}[thm]{Definition}

\newtheorem{cor}[thm]{Corollary}

\newtheorem{lemma}[thm]{Lemma}

\def\ap1{{\square}}

\newcounter{smallromans}

\newcounter{smallarabics}

\newcounter{smallalphs}

\newcommand\bproof{\begin{proof}}
\newcommand\eproof{\end{proof}}
\newcommand{\beginpf}{\smallskip\textbf{Proof. }}
\newcommand{\eopf}{\hfill $\Box$}
\def\cA{{\mathcal A}}
\def\cB{{\mathcal B}}

\def\cal#1{{\mathcal {#1}}}

\begin{document}

\title[The bidual of a radical operator algebra]{The Bidual of a Radical Operator algebra can be Semisimple}
%Operator algebras hereditary in their 
%bidual]{Operator algebras which are hereditary subalgebras of their bidual}
\author{Charles John Read}

\address{Department of Pure Mathematics,
University of Leeds,
Leeds LS2 9JT,
England}
 \email[Charles John Read]{read@maths.leeds.ac.uk}
\thanks{Read is grateful for support from UK research council grant EP/K019546/1, and for helpful suggestions from David Blecher.}

\begin{abstract}
The paper of S. Gulick [Sidney (Denny) L. Gulick, {\it Commutativity and ideals in the biduals of topological algebras}, Pacific J. Math {\bf 18} No. 1, 1966] contains some good mathematics, but it also contains an error. It claims that for a Banach algebra $A$, the intersection of the Jacobson radical of $A^{**}$ with $A$ is precisely the radical of $A$ (this is claimed for either of the Arens products on $A^{**}$ - in itself a reasonable claim, because $A$ is always contained in the topological centre of $A^{**}$, so a fixed $a\in\cA$ lies in the radical of $A^{**}$ with the first Arens product, if and only if it lies in the radical of $A^{**}$ when that Banach space is given the second Arens product, if and only if $ab$ is quasinilpotent for every $b\in A^{**}$). In this paper we begin with a simple counterexample to that claim, in which $A$ is a radical operator algebra, but not every element of $A$ lies in the radical of $A^{**}$. We then develope a more complicated example $\cA$ which, once again, is a radical operator algebra, but $\cA^{**}$ is semisimple. So $\rad \cA^{**}\cap \cA$ is zero, but $\rad \cA=\cA$. We conclude by examining the uses Gulick's paper has been put to since 1966 (at least 8 subsequent papers refer to it), and we find that most authors have used the correct material from that paper, and avoided using the wrong result. We reckon, then, that we are not the first to suspect that the result $\rad A^{**}\cap A=\rad A$ was wrong; but we believe we are the first to provide ``neat'' counterexamples as described.  
 \end{abstract}

\maketitle

\section{Introduction}  
The theorem in which Gulick makes the claim $\rad A^{**}\cap A=\rad A$ is Theorem 4.6 of \cite{G}.
We believe that the place where his proof breaks down is nearby, in the proof of Lemma 4.5, the seventh line: ``note that $M_E$ is once again a maximal regular left ideal in $E$''. We could not see why this should be so, and Theorem 4.6 is definitely false; this introductory section contains a counterexample.

We shall always be working with operator algebras (norm closed subalgebras of the algebra $B(H)$ of all operators on a Hilbert space $H$), so the question of which Arens product is involved need never be addressed, for as is well known, every operator algebra is Arens regular - the two products coincide. 

Let us conclude this Introduction with 
the simpler counterexample mentioned in the Abstract.

Let $H$ be a Hilbert space with orthonormal basis $(e_i)_{i\in\bN}$. Let $T_0:H\to H$ be the operator with 
\beq\label {1}
T_0e_i=\bcas
e_{i+1},&\textrm{if $i$ is odd;}\cr
0,&\textrm{if $i$ is even.}
\ecas
\eeq
For $n\in\bN$, let $T_n:H\to H$ be the rank 1 operator with
 \beq\label{2}
T_ne_i=\bcas
e_{i+1},&\textrm{if $i=2n$;}\cr
0,&\textrm{otherwise.}
\ecas
\eeq
Let $A$ denote the operator algebra (the norm-closed subalgebra of $B(H)$) generated by $\{T_n:n\in\bN_0\}$.
\blem\label{1} $A$ is radical. 
\elem
\proof
First, $T_0^2=0$ and each $T_n\ (n\ge 1)$ has rank 1, so everything in $A$ is of form $\lam T_0+K$, where $\lam\in\bC$ and $K$ is a compact operator. Second, the subspaces $E_k=\linbar\{e_i:i>k\}\subset H$ are invariant for every $T_n$ (and hence for every $T\in A$); indeed, every $T\in A$ maps $E_k$ into $E_{k+1}$ ($k\in\bN_0$). So, let $T=\lam T_0+K\in A$, with $\lam\in\bC$ and $K\in K(H)$. It is enough to show that $T$ is quasinilpotent. Since $K$ is compact, the norms $\veps_n=\nm{K|_{E_n}}$ tend to zero as $n\to\infty$. 
Furthermore, since $T_0^2=0$, we have 
\beq\label{3}
\nm{T^2|_{E_n}}=\nm{\lam T_0K+\lam KT_o+K^2|_{E_n}}\le 2|\lam|\veps_n+\veps_n^2=\del_n,
\eeq
with $\del_n\to 0$ as $n\to\infty$. Now $T^{2k}= T^2|_{E_{2k-2}}T^2|_{E_{2k-4}}\ldots T^2|_{E_2}T^2|_{E_0}$, hence $\nm{T^{2k}}\le \prodd j0{k-1}\del_{2j}$, so $\nm{T^{2k}}^{1/k}\to 0$. Plainly $T^2$, and hence $T$ itself, is quasinilpotent.
\endproof

\bthm\label{t1} $T_0\notin\rad A^{**}$, so $A=\rad A\subsetneq A\cap\rad A^{**}$.
\ethm
\proof Now $A\subset B(H)$, and $B(H)$ is of course a dual Banach algebra, so there is a natural projection 
from $B(H)^{**}$ (the third dual of the Banach space of trace class operators on $H$) onto $B(H)$. This projection is an algebra homomorphism, so when we restrict it to $A^{**}\subset B(H)^{**}$, we get a representation of $A^{**}$ acting on $H$, such that the canonical image $A\subset A^{**}$ acts on $H$ in its usual way, and the representation of $A^{**}$ consists of the weak-* closure of $A$ in $B(H)$. 

Among the operators in this weak-* closure is the weak-* convergent sum $T=\summ n1\infty T_n$, with 
\beq\label{4}
Te_i=\bcas
e_{i+1},&\textrm{if $i$ is even}\cr
0,&\textrm{if $i$ is odd}
\ecas
\eeq 
The product $TT_0$ has $TT_0e_i=e_{i+2}$ (if $i$ is odd) or zero (if $i$ is even); so $\nm{(TT_0)^k}=1$ for all $k$, indeed 1 is in the spectrum of $TT_0$. If $\tau\in A^{**}$ is any element represented as $T$ by this representation, then $1\in$ Sp$(\tau T_0)$. So $T_0$ does not lie in the Jacobson radical of $A^{**}$, by a well known characterization of that radical.  
\endproof

 Note that the proof given above does not depend on the faithfulness (injectivity) of the natural representation of $A^{**}$ in $B(H)$. However, when we give the more complicated  counterexample - when we make the claim that the bidual of our radical algebra $\cA$ is semisimple - we will have to show that the analogous  representation for the bidual of that algebra is indeed faithful.

\section{The main construction}

We now seek to develope the example given in the Introduction, into an example $\cA$ where $\cA$ is radical, but $\cA^{**}$ is semisimple. 

\bdfn\label{d1}
Let $S$ denote the free unital semigroup on two generators $g,h$. If $s\in S$ with $s=\gam_n\gam_{n-1}\ldots \gam_2\gam_1=\prodd j0{n-1}\gam_{n-j}$, and each $\gam_i\in\{g,h\}$, we define the length $l(s)=n$; the depth $\rho(s)=\#\{i:1\le i\le n,\gam_i=h\}$. If $n>0$ (that is, if $s\ne 1$, the unit), the predecessor $p(s)=\prodd j1{n-1}\gam_{n-j}$. We define $S^-=S\setminus\{1\}$. 

We define the Cayley graph $G$ of $S$ to be an abstract directed graph with vertex set $S$, and a directed edge $p(s)\to s$ for each $s\in S^-$.   
\edfn
Note that $G$ is an infinite tree with root vertex 1, such that every vertex $s\in S$ has two outward edges leaving it (the edges $s\to gs$ and $s\to hs$), and every vertex $s\in S^-$ not equal to the root vertex, has a single edge entering it (the edge $p(s)\to s$). If $l(s)=k$, the unique directed path from 1 to $s$ consists of 
$k+1$ vertices $1\to p^{k-1}(s)\to p^{k-2}(s)\to\ldots\to p(s)\to s$. 
\bdfn\label{d2}
For $s\in S$ we define the weight $w(s)=2^{-\rho(s)}$, and if $l(s)=l$ we define 
\beq\label{5}
W(s)=\prodd j0{l-1} w(p^js).
\eeq
We define a Hilbert space $\cH=l^2(S,W)$ to be the collection of all formal sums $\bx=\sum_{s\in S}x_s\cdot s$ with $x_s\in\bC$ ($s\in S$), and 
\beq\label{6}
\nm\bx^2=\sum_{s\in S}W(s)^2|x_s|^2<\infty.
\eeq
We define a particular subset $\cC\subset S^-$, the ``colour set''
\beq\label{7}
\cC=\{g^k:k\in\bN\}\cup\{g^khs:k\in\bN_0, s\in S, 1+l(s)\divides k\}.
\eeq
(here and elsewhere we use ``$1+l(s)\divides k$'' for ``$1+l(s)$ divides $k$'').
We define a ``colour map'' $\mu:S^-\to \cC$ recursively as follows:
\beq\label{8}
\mu(s)=\bcas
s,&\textrm{if $s\in\cC$;}\cr
\mu(p^{n-k'}y),&\textrm{if $s=g^khy$, $y\in S$, $l(y)=n$}, \cr
&\textrm{$1\le k'\le n$, $k\equiv k'$ (mod $n+1$)}.
\ecas
\eeq
\edfn
Note that equation \refeq{8}\ really ``works'' as a recursive definition, because if $s\notin\cC$, 
we necessarily have $s=g^khy$ for some $k\in\bN$ such that $1+l(y)\notdivides k$; so writing
$n=l(y)$, there is a unique $k'\in[1,n]$ such that $k'\equiv k$ (mod $n+1$). The iterated predecessor 
$p^{n-k'}y$ will not be equal to 1 because $k'>0$ and $l(y)=n$, so $\mu(p^{n-k'}y)$
will be (recursively) defined. 
Note that for $s\in S$, the colour $\mu(hs)$ is always equal to $hs$, while the colour $\mu(gs)$ is either $gs$ itself, or one of the iterated predecessors of $gs$. So we never have $\mu(gs)=\mu(hs)$ for any $s\in S$.
\bdfn\label{d3}
For each colour $c\in\cC$, we define a linear map $T_c\in B(\cH)$ by its action on the basis $S$, as follows:
for each $s\in S$, we define
\beq\label{9}
T_c(s)=\bcas gs,&\textrm{if $\mu(gs)=c$;}\cr
 hs,&\textrm{if $\mu(hs)=c$;}\cr
 0,&\textrm{otherwise.}
\ecas
\eeq
\edfn

Each $T_c$ is a weighted shift operator (for $S$ is an orthogonal, though not an orthonormal, basis of $H$).
Writing $e_s=W(s)^{-1}\cdot s$ ($s\in S$) for the corresponding orthonormal basis, and 
 giving due regard to the fact that $W(s)/W(p(s))=w(s)$ for each $s\in S^-$, we have 
\beq\label{10}
T_c(e_s)=\bcas w(gs)e_{gs},&\textrm{if $\mu(gs)=c$;}\cr
 w(hs)e_{hs},&\textrm{if $\mu(hs)=c$;}\cr
 0,&\textrm{otherwise.}
\ecas
\eeq
This implies that for each $c\in\cC$,
\beq\label{11}
\nm{T_c}=\max\{w(x):\mu(x)=c\}=w(c)=2^{-\rho(c)}.
\eeq
\bdfn\label{d4}
We define two families of ``coordinatewise'' orthogonal projections on $\cH$. For $n\in\bN_0$, $P_n$ is the orthogonal projection onto $\linbar\{s\in S:\rho(s)=n\}$, and $\bar P_n=\summ i0n P_i$; while $\pi_n$ is the orthogonal projection onto $\lin\{s\in S:l(s)=n\}$, and $\bar\pi_n=\summ i0n\pi_n$. 

We also define, for $n\in\bN_0$, a subgraph $G\up n$ of $G$, obtained from $G$ by deleting some of the edges. Specifically, $G\up n$ is a graph with vertex set $S$, and a directed edge $p(s)\to s$ for every $s\in S$ such that the ``colour depth'' $\rho\mu(s)\le n$. (Equivalently, we obtain $G\up n$ by deleting from $G$ every edge $p(s)\to s$ such that the colour depth $\rho\mu(s)$ is greater than $n$).
If $K\subset G\up n$ is a connected component, we define the coordinatewise projection $Q_{n,K}$ by 
\beq\label{12}
Q_{n,K}(s)=
\bcas s,&\textrm{if $s\in K$;}\cr
0,&\textrm{otherwise;}
\ecas \hskip 1cm (s\in S).
\eeq
We define $H_{n,K}=Q_{n,K}(\cH)$.
\edfn 
Note that while $\pi_n$ has finite rank $2^n$, the projection $P_n$ always has infinite rank (even when $n=0$, when it is the orthogonal projection onto $\linbar\{g^k:k\ge 0\}$).  

 \bdfn\label{d5} We define an algebra $\cA_0\subset B(\cH)$. $\cA_0$ is the non-unital subalgebra of $B(\cH)$ generated by the operators $T_c$ ($c\in\cC$). We define the operator algebra $\cA=\overline{A_0}$, the norm closure of $A_0$ in $B(\cH)$. We define $\cA\up n\subset \cA_0$ to be the linear span of products $T=T_{c_k}T_{c_{k-1}}\ldots T_{c_2}T_{c_1}=\prodd i0{k-1}T_{c_{k-i}}$ such that $c_i\in\cC$ and $\max\{\rho(c_i): 1\le i\le k\}=n$. We define $\bar\cA\up n=\summ r0n\cA\up r$, the subalgebra of $\cA_0$ generated by maps $T_c$ ($c\in\cC$) with $\rho(c)\le n$. 

For  $n,r\ge 0$, let $S_{n,r}=\{s\in S:$ the path from 1 to $s$ in $G$ contains exactly $r$ edges $p(u)\to u$ with colour depth $\rho\mu(u)>n\}$. Let $P_{n,r}$ be the orthogonal projection onto $\linbar(S_{n,r})$, and let $\bar P_{n,r}=\summ t0rP_{n,t}$.
\edfn

Note that $S_{n,0}=\{s\in S:\rho(s)\le n\}$, so $P_{n,0}=\bar P_n$ for each $n\in\bN_0$.

\blem\label{l2}
\ni(a) For each $n\in\bN_0$, the subspaces $\ker \bar P_n$, $\ker\bar\pi_n\subset \cH$ are invariant for $\cA$. Further, $\cA$ maps $\ker\bar\pi_n$ into $\ker\bar\pi_{n+1}$ for each $n$.

\ni(b) For each component $K$ of $G\up n$, the subspace $H_{n,K}$ is invariant for $\bar\cA\up n$ and also for the hermitian conjugate $(\bar\cA\up n)^*$. The component of $G\up n$ containing 1 is $S_{n,0}$, and the associated projection is $\bar P_n$. 

\ni (c) Every map $T_c$ with $\rho(c)>n$ maps $\cH$ into $\ker\bar P_n$.

\ni (d) For $T\in \cA_0$, the decomposition $T=\summ n1\infty T\up n$, with $T\up n\in \cA\up n$, is unique and continuous; writing $\bar T\up n=\summ i0nT\up i$, we have $\nm{\bar T\up n}\le \nm T$ for every $n$ and $T$; in fact $\bar T\up n=\summ r0\infty P_{n,r}TP_{n,r}$ in the strong operator topology, while $T-\bar T\up n=\summ r0\infty (1-\bar P_{n,r})TP_{n,r}$.

\ni (e) For all $s\in S$ we have $\rho\mu(s)\le \rho(s)$, with equality if $s\in hS$.
\elem
\ni\proof (a) is obvious because the generating maps $T_c$ all map an element $s\in S$ to $gs, hs,$ or zero; and we have $\rho(gs)\ge \rho(s),\ \rho(hs)\ge \rho(s)$, and $l(gs)=l(s)+1$, $l(hs)=l(s)+1$ for all $s\in S$. 

For $c\in\cC$, we have $\sprod{T_cs}t\ne 0$ ($s,t\in S$) only when there is an edge $s\to t$ in $G$, and  $\mu(t)=c$. So if $T\in\bar\cA\up n$, the algebra generated by maps $T_c$ with $\rho(c)\le n$, and if $\sprod {Ts}t\ne 0$, then there is a path from $s$ to $t$ in $G$, and each edge $p(u)\to u$ in that path has $\rho\mu(u)\le n$, so the edge $p(u)\to u$ is present in the graph $G\up n$. Thus $s,t$ belong to the same component of $G\up n$. So for a connected component $K\subset G\up n$, the associated subspace  $H_{n,K}$ is invariant for both $\cA\up n$ and   $(\cA\up n)^*$, establishing the first part of (b).

The component of $G\up n$ containing 1 is the set of $s\in S$ such that the path from 1 to $s$ in $G$ contains only edges $p(u)\to u$ with $\rho\mu(u)\le n$. Now for any $u\in S$, $\mu(u)$ is either $u$ itself or one of the iterated predecessors $p^i(u)$; taking predecessors cannot increase the depth $\rho(u)$, so $\rho\mu(u)
\le\rho(u)$ for all $u$. If $s\in S$ with $\rho(s)\le n$, then every edge $p(u)\to u$ in the path from 1 to $s$ has colour depth $\rho\mu(u)\le n$ also, so $s$ lies in the component of $G\up n$ containing 1. Conversely, if $\rho(s)>n$ then we have $s=g^kht$ for some $t\in S$ and $k\in \bN_0$; the edge $t\to ht$ is part of the path from 1 to $s$, and 
$ht\in\cC$ by \refeq{7}, so the colour depth $\rho\mu(ht)=\rho(ht)=\rho(s)>n$, therefore 
$s$ is not in the connected component of $G\up n$ containing 1. Therefore that component is precisely $\{s:\rho(s)\le n\}$, and the associated coordinatewise projection is $\bar P_n$. Thus we have established the second part of (b), and also part (e).

For part (c), note that $T_c$ maps $\cH$ into $\linbar\{x\in S:\mu(x)=c\}$; if $\rho(c)>n$ then this subspace is contained in $\linbar\{x\in S:\rho\mu(x)>n\}\subset\linbar\{x\in S:\rho(x)>n\}$ (by part (e)),
$\subset\ker\bar P_n$. 

To prove part (d), we note that the edges of $G\up n$ include the edge $p(u)\to u$ only if $\rho\mu(u)\le n$, hence the set $S_{n,r}$ is a union of some of the components $K$ of $G\up n$. So by part (b) of this Lemma, each image $P_{n,r}\cH$ is $\cA\up n$ invariant; but for $c\in\cC$ with $\rho(c)>n$, $T_c$ maps $P_{n,r}\cH$ into $P_{n,r+1}\cH$ because $\sprod {T_cs}{t}\ne 0$ ($s,t\in S$) only when $s=p(t)$ and the colour depth $\rho\mu(t)>n$. Now take any $T\in \cA_0$ and write $T=\sum_iT\up i$ with each $T\up i\in\cA\up i$.  We have $\bar T\up n=\summ i0nT\up i\in\bar\cA\up n$, so each $P_{n,r}\cH$ is a $\bar T\up n$-invariant subspace; but $T-\bar T\up n$ maps  $P_{n,r}\cH$ into $\opluss i{r+1}\infty P_{n,i}\cH$. Therefore we have 
\beq\label{13}
\bar T\up n=\summ r0\infty \bar T\up nP_{n,r}=\summ r0\infty  P_{n,r}T\up nP_{n,r}=\summ r0\infty P_{n,r}TP_{n,r},
\eeq
while $T-\bar T\up n=\summ r1\infty (1-\bar P_{n,r})TP_{n,r}$
as required by the Lemma.
This shows that the decomposition $T=\summ i0\infty T\up i$ is indeed unique, and furthermore the compression $\bar T\up n$ as given by \refeq{13}\ plainly satisfies $\nm{\bar T\up n}\le \nm T$. 
Thus the lemma is proved. \endproof

\bdfn\label{d3b}
Let us write $\cB\up n$ ($\bar\cB\up n$) for the norm closure of $\cA\up n$ ($\bar\cA\up n$) in $B(\cH)$. Let us write $\Del_n$ for the map $\cA_0\to\cA\up n$ with
$\Del_n(T)$ the unique element $T\up n\in\cA\up n$ such that $T=\summ n0\infty T\up n$; and let
 $\bar\Del_n:\cA_0\to\bar\cA\up n$ be the map  $\summ i0n\Del_i$.
\edfn

The maps $\Del_n,\bar\Del_n$ are uniformly norm bounded by part (d) of the previous lemma; so they extend continuously to maps $\Del_n:\cA\to\cB\up n$ and $\bar\Del_n:\cA\to\bar\cB\up n$; and because of the uniform bound on $\nm{\bar\Del_n}$ (each $\bar\Del_n$ is contractive), we have 
$T=\summ n0\infty\Del_nT=\summ n0\infty T\up n$, with $T\up n\in\cB\up n$, for all $T\in\cA$. The formulae $\bar\Del_nT=\bar T\up n=\summ r0\infty P_{n,r}TP_{n,r}$ and $T-\bar T\up n=\summ r0\infty (1-\bar P_{n,r})TP_{n,r}$ remain true in the strong operator topology.

\section{$\cA$ is radical.}
  In order to prove that our algebra $\cA$ is radical, the main theorem we need is the following:
\bthm\label{t2}
Every $T\in\cA\up n$, or the norm closure thereof, satisfies
\beq\label{14}
(1-\bar\pi_k)\bar P_nT\to 0 \textrm{ as } k\to\infty.
\eeq
Indeed, $\bar P_nT$ is a compact operator.
Furthermore, every $T\in\bar\cA\up n$ satisfies
\beq\label{15}
\nm T=\nm{\bar P_nT\bar P_n}.
\eeq
\ethm

Let us prove the first part of the Theorem. From \refdfn{d2}, we find that if $s\ne c$ but $\mu(s)=c$, then we  must have $s=g^khyc$ for some $k\in\bN_0$ and $y\in S$. In particular, $\rho(s)>\rho(c)$. 
So if $\rho(c)=n$, then the map $\bar P_nT_c$ in fact has rank 1; it maps $p(c)$ to $c$, and all other $s\in S$ to zero. Any product $T=\prodd i0{k-1} T_{c_{k-i}}$ with $c_i\in\cC$ and $\max\{\rho(c_i):1\le i\le k\}=n$, accordingly satisfies $\bar P_nT=\prodd i0{k-1} \bar P_nT_{c_{k-i}}$ (because $\ker\bar P_n$ is an invariant subspace for each $T_{c_j}$), so the rank of $\bar P_nT$ is at most 1. 
$\cA\up n$ is the linear span of such maps, so any $T\in\cA\up n$, or its norm closure, will have 
$\bar P_nT$ a compact operator; hence, $ \nm{(1-\bar\pi_k)\bar P_nT}\to 0$ as $k\to\infty$.\endproof

To prove the second part of the Theorem, we  need certain preliminaries, which we bring together in the following Lemma:

\blem\label{l3}
\ni(a) Let $K$ be a connected component of $G\up n$. Then either   $K=S_{n,0}$, the component which contains 1, or $K$ consists of a path $y\to gy\to g^2y\to\ldots g^my$ for some $y\in S$ and $m\in\bN$ such that the colour depths $\rho\mu(y)>n$, $\rho\mu(g^{m+1}y)>n$, but $\rho\mu(g^iy)\le n$ for $i\in [1,m]$. Furthermore, there is a path $s_0\to s_1\to \ldots\to s_m$ in the component  $S_{n,0}$ such that the colours $\mu(s_i)=\mu(g^iy)$ for each $i\in[1,m]$. 

\ni(b) Let $M\in M_{m+1}(\bC)$ be a strictly lower triangular matrix, and let $\Nm$ and $\Nm'$ be two norms on $\bC^{m+1}$, with $\nm{\lam_0,\lam_1,\ldots \lam_m}=(\summ i0m\omega_i^2|\lam_i|^2)^{1/2}$ and $\nm{\lam_0,\lam_1,\ldots \lam_m}'=(\summ i0m(\omega_i')^2|\lam_i|^2)^{1/2}$ for positive constants $\omega_i,\omega_i'$ ($\range 1im$). 
Suppose we have 
\beq\label{16}
\frac{\omega_{i+1}'}{\omega_i'}\le\frac 12\cdot\frac{\omega_{i+1}}{\omega_i}
\eeq
for each $\range i0{m-1}$. Then 
\beq\label{17}
\nm M'\le \nm M.
\eeq

\ni(c) For every $T\in\cA$, we have $\nm{(1-\bar P_n)T}\to 0$ as $n\to\infty$.
\elem

\ni{\bf Proof of Lemma:} (a) Suppose $K\ne S_{n,0}$. Since $K$ cannot meet $S_{n,0}$, every vertex $x\in K$ must have $\rho(x)>n$. But if $x\to x'$ is an edge in $K$, we must have $\rho\mu(x')\le n$, therefore $\mu(x')\ne x'$, so $x'\notin\cC$, so $x'=g^khz$ for some $z\in S$ and $k>0$ with $\rho(hz)=\rho(x')>n$. Indeed, we must have $1+l(z)\notdivides k$. Every edge of $K$ must be of form $x\to gx$ rather than $x\to hx$, so $K$ does indeed consist of a path (finite or infinite) of form $g^rhz\to g^{r+1}hz\to g^{r+2}hz\ldots$, for some $r\ge 0$. But we have the condition  $1+l(z)\notdivides k$ for any $k$ such that $k>r$ and $g^khz$ is in the path; so the path is finite. Its last vertex must be $g^thz$ for some $t$ with $t-r\le 1+l(z)$.
Writing $m=t-r$ and $y=g^rhz$ we see that $K=\{g^iy:\range i0m\}$.

If $r>0$, we must have $\rho\mu (y)=\rho\mu(g^rhz)>n$ or we could continue the path in $K$ backwards to include the vertex $g^{r-1}hy$. If $r=0$, we have $\mu(y)=\mu(hz)=hz$ so $\rho\mu(y)>n$ anyway. Also,  we must have   $\rho\mu(g^{t+1}hy)>n$ or we could include the vertex $ g^{t+1}hy$ in our component $K$. For $i\in(r,t]$ we have    $\rho\mu(g^ihy)\le n$ because the edge $
g^{i-1}hy\to g^ihy$ lies in $K$. Thus the component $K$ is as described in part (a) of this Lemma.

 To complete the proof of part (a), we claim that there is a sequence $s_0\to s_1\to\ldots s_m\in S_{n,0}$ such that $\mu(s_i)=\mu(g^iy)$ for each $i\in[1,m]$. This is proved by induction on $l(y)=\min\{l(u):u\in K\}$. If $l(y)\le n$, there is nothing to prove because the component is $S_{n,0}$ after all. If the component $K\ne S_{n,0}$, write $K=\{g^ihz:r\le i\le r+m\}$.
We return to equation \refeq{8} to compute the colours $\mu(g^ihz)$ for $i\in (r,r+m]$. Writing $l=l(z)$ and  $z=\prodd i0{l-1} z_{l-i}$ ($z_j\in\{g,h\}$),  we find that if $i'\in[1,l]$ is the unique integer with $i'\equiv i$ (mod $1+l$), then
$\mu(g^ihz)=\mu(\prodd j0{i'-1}z_{i'-j})=\mu(p^{l-i'}z)$. If $r_0\in[0,l]$ satisfies $r_0\equiv r$ (mod $l+1$), then the sequence  $\mu(g^iy)$ ($\range i1m$) is the sequence $\mu(p^{l-r_0-i}z)$ ($\range i1m$). The vertices $ (p^{l-r_0-i}z)_{i=0}^m$ form a path in $G$ which, since it involves the same colours for $i>0$, is also a path in $G\up n$. So this path is part of a component $K'$ of $G\up n$. If $K'=S_{n,0}$ we are done; if not, we note that the minimum length of an element of $K'$ is strictly less than $l(y)$, so the result follows by induction hypothesis.

(b)  Let $(e_i)_{i=0}^m$ be the unit vectors of $\bC^{m+1}$, and write 
$Me_i=\sum_{j>i}M_{j,i}e_j$.  We may assume $\nm M=1$, in which case 
$|M_{j,i}|\le \nm {e_i}/\nm{e_j}=\omega_i/\omega_j$ for all $i$ and $j$.
For $k\in[1,m]$, the weighted shift matrix $M\up k$ with 
\beq\label{18}
M\up ke_i=\bcas
M_{i+k,i}e_{i+k},&\textrm{ if } i+k\le m;\cr
0,&\textrm{ if } i+k>m;
\ecas
\eeq    
satisfies $\nm{M\up k}'=\max_{i\in[0,m-k]}|M_{i+k,i}|\omega_{i+k}'/\omega_i'$
$\le \max_{i\in[0,m-k]}(\omega_i/\omega_{i+k})\cdot \omega_{i+k}'/\omega_i'$ 
$\le 2^{-k}$ by \refeq{16}. Then $M=\summ k1mM\up k$ so $\nm M'\le\summ k1m2^{-k}< 1$. 

(c) Let $c\in\cC$. From \refeq{10}, for $s\in S$ we have 
$$
(1-\bar P_n)T_ce_s=\bcas
w(gs)e_{gs},&\textrm{if $\mu(gs)=c$ and $\rho(gs)>n$}\cr
w(hs)e_{hs},&\textrm{if $\mu(hs)=c$ and $\rho(hs)>n$}\cr
0,&\textrm{otherwise}.
\ecas
$$
But $w(x)=2^{-\rho(x)}$, so $\nm{(1-\bar P_n)T_c}\le 2^{-n-1}$. We will also have $\nm{(1-\bar P_n)T}\to 0$ for any operator $T$ in the norm closed right ideal generated by the operators $T_c$. But this right ideal is the entire algebra $\cA$. 
\endproof

{\bf Proof of \refthm{t2}, second part:} By \reflem{l2}\ part (b), when $T\in\bar \cA\up n$ we have
$T=\sum_KQ_{n,K}TQ_{n,K}$, where the sum is taken over the connected components $K$ of $G\up n$. So, 
\beq\label{19}
\nm T=\sup_K\nm{Q_{n,K}TQ_{n,K}}.
\eeq
 If $K=S_{n,0}$, the component containing 1, then the norm $\nm{Q_{n,K}TQ_{n,K}}=\nm{\bar P_nT\bar P_n}$. If $K$ is any other component, we claim that the norm is at most $\nm{\bar P_nT\bar P_n}$. By \reflem{l3} part (a), we can write 
$K=\{g^iy:0\le i\le m\}$ for suitable $y\in S$ and $m$; writing $\gam_i$ for the colour $\mu(g^iy)$, there is also a set $\kappa=\{s_i:0\le i\le m\}\subset S_{n,0}$ such that the colour $\mu(s_i)=\gam_i$ for $i\in[1,m]$. Let $q$ denote the orthogonal (coordinatewise) projection onto $\lin(\kappa)$. 
If $c_1,c_2,\ldots c_r\in \cC$, then the compression $\tau_1=Q_{n,K}T_{c_r}T_{c_{r-1}}\ldots T_{c_1}Q_{n,K}$ sends $g^iy$ to $g^{i+r}y$, if $i+r\le m$ and $c_i=\gam_{r+i}$ for each $\range i1r$; otherwise, we have   $\tau_1g^iy=0$. Similarly, the compression $\tau_2=q T_{c_r}T_{c_{r-1}}\ldots T_{c_1}q$ sends $s_i$ to $s_{i+r}$ if $i+r\le m$ and 
$c_i=\gam_{r+i}$ for each $\range i1r$; otherwise, we have  $\tau_2s_i=0$. So the compressions $\tau_1$ and $\tau_2$ are intertwined by the map $\eta$ sending $g^iy$ to $s_i$ for each $i$. Indeed, if $T\in\bar\cA\up n$, the compressions $\tau= Q_{n,K}TQ_{n,K}$ and $\tau'=qTq$ are intertwined, with $\eta\tau=\tau'\eta$. So $\tau$ has the same $(m+1)\times(m+1)$ matrix $M$ with respect to the basis $(g^iy)_{i=0}^m$ of $Q_{n,K}\cH$, as $\tau'$ has with respect to the basis $(s_i)_{i=0}^m$ of $q\cH$. $M$ is strictly lower triangular, because all such compressions $qTq$ map $s_i$ into $\lin\{s_j:j>i\}$ for each $i$. 
The norm on $q\cH$ is given by 
$\nm{\summ i0m\lam_is_i}=(\summ i0m\omega_i^2|\lam_i|^2)^{1/2}$, where $\omega_i=W(s_i)$. 
 The norm on $Q_{n,K}\cH$ is likewise given by 
$\nm{\summ i0m\lam_ig^iy}=(\summ i0m(\omega_i')^2 |\lam_i|^2)^{1/2}$, where $\omega_i'=W(g^iy)$. For $0\le i<m$, the ratio $\omega_{i+1}/\omega_i= W(s_{i+1})/W(s_i)=w(s_{i+1})$ because there is an edge $s_i\to s_{i+1}$ in $G$; and $w(s_{i+1})\ge 2^{-n}$ because $s_{i+1}\in S_{n,0}$ so $\rho(s_{i+1})\le n$. 
On the other hand, the ratio   $\omega_{i+1}'/\omega_i'= W(g^{i+1}y)/W(g^iy)=w(g^{i+1}y)\le 2^{-n-1}$, because $g^{i+1}y\notin S_{n,0}$ so $\rho(g^{i+1}y)\ge n+1$. We deduce that 
$\omega_{i+1}'/\omega_i'\le\frac 12\cdot \omega_{i+1}/\omega_i$. By \reflem{l3} part (b), we have 
$\nm{Q_{n,K}TQ_{n,K}}=\nm{\tau}\le \nm{\tau'}$, and of course $\nm{\tau'}\le\nm{\bar P_nT\bar P_n}$ because the orthogonal projection $q\le \bar P_n$.  By \refeq{19}, the norm of $T$ is the supremum of $\nm{\bar P_nT\bar P_n}$ and the norms $\nm{Q_{n,K}TQ_{n,K}}$ for all other connected components $K\subset G\up n$; so $\nm T=\nm {\bar P_nT\bar P_n}$ as claimed by the Theorem. \endproof

We can now prove the main theorem of this section:

\bthm\label{t3}
$\cA$ is radical.
\ethm
\proof If not, let $T\in\cA$ have spectral radius at least 1. By \reflem{l3}, there is an $n\in\bN$ such that $\nm{(1-\bar P_n)T}\le \frac 12$. I claim that the spectral radius of the compression $\bar P_nT\bar P_n$ is at least 1. For by \reflem{l2}(a), for each $k\in\bN$ we have  $T^k=\bar P_nT^k\bar P_n+(1-\bar P_n)T^k$ (any $k\in\bN$) because $\ker\bar P_n$ is an invariant subspace for $\cA$; indeed,
$T^k=(\bar P_nT\bar P_n)^k+(1-\bar P_n)T^k$, because the compression map $T\to \bar P_nT\bar P_n$ is an algebra homomorphism on $\cA$. So for all $k>0$, $T^k=(\bar P_nT\bar P_n)^k+(1-\bar P_n)T\cdot T^{k-1}$,
hence
$$
\nm{T^k}\le \nm{(\bar P_nT\bar P_n)^k}+\frac 12\cdot\nm{T^{k-1}}
\le \nm{(\bar P_nT\bar P_n)^k}+\frac 12\cdot\nm{(\bar P_nT\bar P_n)^{k-1}}+\frac 14\cdot\nm{T^{k-2}}
$$
$$
\le\ldots\le 2^{-k}+\summ j0{k-1} 2^{-j}\nm{(\bar P_nT\bar P_n)^{k-j}}.
$$
If the spectral radius of $\bar P_nT\bar P_n$ is less than 1, we can find $r<1$ and $C>0$ such that 
$\nm{(\bar P_nT\bar P_n)^{j}}\le Cr^j$ for all $j\in\bN$, so we have 
$1\le \nm{T^k}\le 2^{-k}+\summ  j0{k-1} C\cdot 2^{-j}\cdot r^{k-j}$ $\le 2^{-k}+kC\max(\frac 12, r)^k$ for all $k\in\bN$. This is a contradiction for large $k$, so the spectral radius of the compression $\bar P_nT\bar P_n$ must be at least 1.  

It is thus sufficient to show that for each $T\in\cA$ and $n\in\bN$, the compression $\bar P_nT\bar P_n$ is quasinilpotent. Let us prove this by induction on $n$, beginning with the not-quite-trivial case $n=0$.

By \reflem{l2}(d) (and its generalization to $T\in\cA$ rather than $T\in\cA_0$  as discussed after \refdfn{d3b}), we have $\bar P_0T\bar P_0=\bar P_0\bar T\up 0\bar P_0=\bar P_0T\up 0\bar P_0$ for any $T\in\cA$; and $T\up 0\in\cB\up 0$. By \refthm{t2}, we have $(1-\bar\pi_k)\bar P_0T\up 0\to 0$, and by \reflem{l2}(a), $T\up 0$ maps $\ker\bar\pi_k$ into $\ker\bar\pi_{k+1}$ for every $k$. Writing $\veps_k=\nm{(1-\bar\pi_k)\bar P_0T\up 0}$, we have $\veps_k\to 0$, and $(\bar P_0T\bar P_0)^k=$
$
(\bar P_0T\up 0\bar P_0)^k=$ $(1-\bar\pi_{k-1})\bar P_0T\up 0(1-\bar\pi_{k-2})\bar P_0T\up 0\cdot(1-\bar\pi_{k-3})\ldots \bar P_0T\up 0(1-\bar\pi_{0})\bar P_0T\up 0\bar P_0,
$
so $\nm{(\bar P_0T\bar P_0)^k} \le \prodd j0{k-1}\veps_j$, hence $\bar P_0T\bar P_0$ is indeed quasinilpotent. 

Proceeding to the case of a general $n\in\bN$, we note that for $T\in\cA$, $\bar P_nT\bar P_n=\bar P_n\bar T\up n\bar P_n=\bar P_n(T\up n+\bar T\up {n-1})\bar P_n$, where $T\up n\in\cB\up n$ and 
$\bar T\up {n-1}\in\bar\cB\up{n-1}$. 

Writing $\tau=\bar T\up {n-1}$, we have $\tau^k\in\bar\cB\up{n-1}$ for all $k$, so by \refthm{t2}, $\nm{\tau^k}=\nm{\bar P_{n-1}\tau^k\bar P_{n-1}}$ for all $k$. But $\ker\bar P_{n-1}$ is an invariant subspace for $\cA$, so 
$\bar P_{n-1}\tau^k\bar P_{n-1}$ $=$ $(\bar P_{n-1}\tau\bar P_{n-1})^k$; and our induction hypothesis tells us that $\bar P_{n-1}\tau\bar P_{n-1}$ is quasinilpotent. So $\nm{\tau^k}^{1/k}\to 0$ as $k\to\infty$, and also $\nm{(\bar P_n\tau\bar P_n)^k}^{1/k}=$ $\nm{\bar P_n\tau^k\bar P_n}^{1/k}\to 0$ as $k\to\infty$. So $\bar P_n\bar T\up {n-1}$ is quasinilpotent. 

Meanwhile $\sig=\bar P_nT\up n$ is a compact operator by \refthm{t2}, satisfying 
$\veps_k=\nm{(1-\bar\pi_k)\sig}\to 0$ as $k\to\infty$; and both $\sig$ and $\tau$ map $\ker\bar\pi_k$ into $\ker\bar\pi_{k+1}$ for each $k$. 

Let us pick an arbitrary $\del>0$ and choose $C>0$ such that $\nm{(\bar P_n\bar T\up {n-1})^k}\le C\cdot \del^k$ for all $k\in\bN_0$. Then for any $k\in\bN$, we have $(\bar P_nT\bar P_n)^k=(\bar P_nT\up {n-1}+\sig)^k\bar P_n$
$$
=\summ r0k\sum_{{i_0+i_1+\ldots i_r=k-r}\atop{i_j\in\bN_0}}(\bar P_nT\up {n-1})^{i_0}
\cdot\prodd j1r\sig\cdot (\bar P_nT\up {n-1})^{i_j}\cdot \bar P_n
$$
and writing $u_j=\summ tjr(1+i_t)-1$, the product from $j=1$ to $r$ is equal to
$\prodd j1r(1-\bar\pi_{u_j})\sig(\bar P_nT\up {n-1})^{i_j}$; so 
$$
\nm{(\bar P_nT\bar P_n)^k}\le\summ r0k\sum_{{i_0+i_1+\ldots i_r=k-r}\atop{i_j\in\bN_0}}
C^{r+1}\del^{k-r}\cdot\prodd j1r\veps_{u_j}.
$$
Now $u_j\ge j-1$ in all cases, so writing $\eta_j=\prodd j1r\veps_{j-1}$, we have 
$$
\nm{(\bar P_nT\bar P_n)^k}\le\summ r0k\sum_{{i_0+i_1+\ldots i_r=k-r}\atop{i_j\in\bN_0}}
C^{r+1}\del^{k-r}\eta_{r}=\summ r0k {k\choose r} C^{r+1}\del^{k-r}\eta_{r}.
$$
But $\eta_r^{1/r}\to 0$, so we can choose $D>0$ such that $\eta_r\le D\cdot (\del/C)^r$ for all $r$; substituting this in the previous equation, we find that  
$\nm{(\bar P_nT\bar P_n)^k}\le\summ r0k {k\choose r} CD\del^{k} = CD\cdot (2\del)^k$. So the spectral radius of $\bar P_nT\bar P_n$ is at most $2\del$; but $\del>0$ was arbitrary, so $\bar P_nT\bar P_n$ is quasinilpotent. Therefore every $T\in\cA$ is quasinilpotent; $\cA$ is a radical Banach algebra.\endproof

\section{$\abws$ is semisimple.}
We wish to prove the second half of our main result, namely that the bidual $\cA^{**}$ is semisimple. We shall do this by showing that 
the weak-* closure $\abws$ of $\cA$ in $B(\cH)$ is semisimple, and then show that the natural representation 
$\theta:\cA^{**}\to B(\cH)$, whose image is $\abws$, is faithful, so that $\cA^{**}$ itself is semisimple. (Our ``natural representation'' is the restriction to $\cA^{**}$ of the natural projection $\cT^{***}\to \cT^*$, where $\cT$ are the trace-class operators on $\cH$, and $\cT^*=B(\cH)$, 
$\cT^{***}=B(\cH)^{**}$). 

In this section, we show that $\abws$, very unlike $\cA$ itself, is semisimple.

\bdfn\label{d4b} Let $\cC^{<\infty}$ denote the collection of all finite sequences $(c_1,c_2,...,c_m)$ of colours $c_i\in\cC$, for $m\in\bN$ (we exclude $m=0$). For $\bc=(c_1,c_2,...,c_m)\in \cC^{<\infty}$, let $T_\bc$ denote the operator $\prodd i1mT_{c_i}\in\cA_0$. Let $S_\cA\subset \cC^{<\infty}$  be the set of $\bc\in\cC$ such that $T_\bc\ne 0$. 
\edfn 
We think of $S_\cA$ as the ``support'' of $\cA$, because 
clearly every $T\in\cA_0$ is equal to a sum 
\beq\label{20}
T=\sum_{\bc\in S_\cA}\lam_\bc \cdot T_\bc,
\eeq
the coefficients $\lam_\bc\in\bC$ being finitely nonzero. 
\blem\label{l4}
Given $T\in\cA_0$, the coefficients $\lam_\bc(T)$ such that $T=\sum_{\bc\in S_\cA}\lam_\bc(T)\cdot T_\bc$ are unique, and they are weak-* continuous linear functionals of $T$. 
\elem
\proof
For $c\in\cC$, equation \refeq{10}\ tells us that $\sprod {T_ce_s}{e_t}\ne 0$ if and only if $s=p(t)$ and the colour $\mu(t)=c$, in which case it is equal to $w(t)$. Any easy induction then tells us that for $\bc=(c_1,c_2,...,c_m)\in S_\cA$,   $\sprod {T_\bc e_s}{e_t}\ne 0$ if and only if $s=p^m(t)$ and, for each $\range i1m$, the colour $\mu(p^{i-1}t)=c_i$. In that case, $\sprod {T_\bc e_s}{e_t}=\prodd i0{m-1}w(p^it)=W(t)/W(s)$. So for fixed $s,t$, the colour sequence $\bc\in S_\cA$ such that $\sprod {T_\bc e_s}{e_t}\ne 0$ is unique if it exists; and since $T_\bc\ne 0$ for $\bc\in S_\cA$, for fixed $\bc\in S_\cA$ there is at least one pair $s,t\in S$ such that $\sprod{T_\bc e_s}{e_t}\ne 0$.

Given $T\in\cA_0$, $T=\sum_{\bc\in S_\cA}\lam_\bc\cdot T_\bc$, we therefore have 
\beq\label{21}
\lam_\bc=\lam_\bc(T)=\frac{W(s)}{W(t)}\sprod {Te_s}{e_t}, 
\eeq
where $s,t$ is any pair such that $\sprod{T_\bc e_s}{e_t}\ne 0$. $\lam_\bc$ is indeed uniquely determined by $T$, and it is indeed a weak-* continuous function of $T$; equation \refeq{21}\ even equates $\lam_\bc\in B(H)_*$ with an element of $\cT$ of rank 1. 
\endproof

Given two elements $\bc=(c_1,\ldots ,c_m),\bd=(d_1,\ldots d_n)$ in $\cC$, we can define the product
$\bc\cdot\bd$ to be the sequence $(c_1,\ldots,c_m,d_1,\ldots,d_n)$. From the equation \refeq{20}, we see that for $T,T'\in\cA_0$, we have 
\beq\label{22}
\lam_\bc(TT')=\sum_{\bd,\be\in S_\cA, \bd\odot\be=\bc}\lam_\bd(T)\cdot \lam_\be(T'),
\eeq
where the product $\bd\odot\be$  denotes concatenation of sequences.
The sum is always finite (it has $m-1$ terms when $\bc=(c_1,\ldots,c_m)$), so equation 
\refeq{22}\ remains true even when we extend $\lam_\bc$ to the weak-* closure $\abws$ of $\cA_0$.

Now for each $c\in\cC$, \refeq{10} tells us that the left support projection $l(T_c)$ for the operator $T_c$ is the orthogonal projection onto $\linbar\{e_t:t\in S^-, \mu(t)=c\}$. 
We also have $\nm {T_c}=w(c)=2^{-\rho(c)}\le 1$.
These left support projections are mutually orthogonal for different colours $c$. The corresponding right support projection $r(T_c)$ is the     projection onto $\linbar\{e_s:s\in S, s=p(t), \mu(t)=c\}$. These right support projections are not mutually orthogonal, but nevertheless, for each $s\in S$ there are only two $t\in S^-$ such that $s=p(t)$, so the norm of any sum $\sum_{c\in S_\cA}\lam_c r(T_c)$ is at most $2\cdot\sup\{|\lam_c|:c\in S_\cA\}$. Hence for any sequence $\bx\in l^\infty(S_\cA)$, the formal  sum 
\beq\label{23}
T=\sum_{c\in S_\cA}\frac {x_c}{w_c}\cdot T_c
\eeq
 satisfies 
$$
T^*T=\sum_{c,d\in S_\cA}  \frac {x_c^*x_d}{w_cw_d}T_c^*l(T_c)l(T_d)T_d
=\sum_{c\in S_\cA}   \frac {|x_c|^2}{w_c^2}T_c^*T_c\le \sum_{c\in S_\cA}  {|x_c|^2}r(T_c),
$$
in particular $\nm{T^*T}\le 2\cdot\nm{\bx}_\infty^2$.
So the sum $T$ in fact converges in the weak-* topology to an element of $\abws$ of norm at most 
$\sqrt 2\cdot\nm{\bx}_{\infty}$.   

\bthm\label{t4} $\abws$ is semisimple.
\ethm
\proof
Let $T\in\abws,\ T\ne 0$. We claim that $T\notin\rad\abws$. Let us choose $s,t\in S$ such that 
$\sprod {Te_s}{e_t}\ne 0$.

Suppose first that $s\ne 1$. Let $l_0=l(s)>0$, and for $\range i1{l_0}$, write $d_i=\mu(p^{i-1}s)=\mu(p^{i+m-1}t)$. Writing $\bd=(d_1,\ldots,d_{l_0})\in \cC^{<\infty}$, we will 
have  $T_\bd(1)=s$ so $\bd\in S_\cA$
and the product $T'=T\cdot T_\bd$ satisfies $\sprod {T'e_1}{e_t}\ne 0$. Furthermore, in order to show 
$T\notin \rad\abws$ it is enough to show that $T'\notin\rad \abws$, because the radical is an ideal. 
So, we can replace $T$ with $T'$ if necessary, and assume that  $\sprod {Te_1}{e_t}\ne 0$.

 Then $\lam_\bc(T)\ne 0$, where $\bc=(c_1,c_2,\ldots c_l)\in S_\cA$ is the unique sequence such that $l=l(t)$, (so $1=p^lt$),  and the colours $\mu(p^{i-1}t)$ ($\range i1l$) are 
$c_i$. Write $\xi_m=g^{(m-1)(l+1)}ht$, and let $E\subset \cC$ be the collection $\{\xi_m:m\in\bN_0\}$ (noting from \refeq{7}\ that these elements are truly elements of the colour set $\cC$). Let us also note that the weight $w_{\xi_m}=2^{-\rho(\xi_m)}=2^{-(1+\rho(t))}$ is independent of $m$.   So $U=\sum_{c\in E}T_c\in \abws$ (for $U$ is a weak-* convergent sum like $T$ in \refeq{23}). We claim that the product $U\cdot T\in\abws$ is not quasinilpotent, so $UT$, and  $T$ itself, are not in the radical of $\abws$.
To prove this, we   compute the inner product $\sprod {(UT)^me_1}{e_{\xi_m}}$ for every $m\in\bN$. Obviously $\lam_\bd(U)=1$ (if $\bd\in E$) or zero otherwise.

Now the length $L=l(\xi_m)=m(1+l)$, and  the colour sequence $\mu(p^{i-1}\xi_m)$ ($\range i1L$)
is obtained from \refeq{8}\ as follows: if $1+l\divides i-1$, we have $p^{i-1}\xi_m=$ $g^{(m-1-r)(l+1)}ht$
$(r=(i-1)/(1+l)\in [0,m))$, and $\mu(p^{i-1}\xi_m)=p^{i-1}\xi_m=\xi_{m-r}\in E$. But if $1+l\notdivides i-1$, then 
writing $i-1=r(l+1)+j$ $(r\in [0,m), j\in [1,l])$, if $r=m-1$ we have 
$p^{i-1}\xi_m=p^jht=p^{j-1}t$, so $\mu(p^{j-1}\xi_m)=c_j$; but if $r<m-1$ we have 
$p^{i-1}\xi_m=g^{(m-2-r)(l+1)+l+1-j}ht$ and the recursive definition in equation \refeq{8}\ tells us that 
$\mu(p^{i-1}\xi_m)=\mu(p^{l-(l+1-j)}t)=\mu(p^{j-1}t)=c_j$ also. So for all $i$, $1\le i\le L$, we have  
\beq\label{24}
\mu(p^{i-1}\xi_m)=\bcas
\xi_{m-r}\in E,&\textrm{ if $i-1=r(1+l)$;}\cr
c_j\notin E,&\textrm{ if $i-1\equiv j$ (mod $l+1$), $1\le j\le l$.} 
\ecas
\eeq
The full sequence $(\mu(p^{i-1}\xi_m))_{i=1}^{L}\in S_\cA$ is the concatenation $\odot_{r=0}^{m-1} (\xi_{m-r}\odot \bc)$, where we slightly abuse notation by writing $\xi_{m-r}$ for the sequence of length 1 in $\cC^{<\infty}$. 
Now from \refeq{21}, the inner product 
\beq\label{25}
\sprod {(UT)^me_1}{e_{\xi_m}}=W(\xi_m)\cdot\lam_{\xi_m}((UT)^m); 
\eeq
and using \refeq{22}\ $2m$ times, we have 
$$
\lam_{\xi_m}((UT)^m)=\sum_{{\bd\up 1,\bc\up 1\ldots \bd\up m,\bc\up m\in S_\cA,}\atop  {\odott i1m (\bd\up i\odot\bc\up i)=\odott r0{m-1} (\xi_{m-r}\odot \bc)}}\prodd i1m \lam_{\bd\up i}(U)\lam_{\bc\up i}(T).
$$
But the coefficient $\lam_{\bd}(U)$ can only be nonzero if the sequence $\bd$ has length 1, and consists of one of the colours $\xi_j\in E$ (in which case the coefficient is equal to 1). There are only $m$ such colours in the sequence 
$\odott r1m (\xi_{m-r}\odot \bc)$, and the rest of the sequence consists precisely of $m$ copies of $\bc$, so in fact
\beq \label{26}
\lam_{\xi_m}((UT)^m)=\lam_\bc(T)^m.
\eeq
Equation \refeq{26}\ makes the rest of the proof rather straightforward. Substituting it in \refeq{25}, we have  $\nm {(UT)^m}\ge |\sprod {(UT)^me_1}{e_{\xi_m}}|=|\lam_\bc(T)|^m\cdot W(\xi_m)$;
where writing $L=m(1+l)$ as usual, we have $W(\xi_m)=\prodd j1Lw(p^{j-1}\xi_m)=2^{-\summ j1L\rho(p^{j-1}\xi_m)}$, from \refdfn{d2}. But $\xi_m=g^{(m-1)(l+1)}ht$ so $\rho(\xi_m)=1+\rho(t)$. And $\rho(p^i\xi_m)\le \rho(\xi_m)$ for all $i\ge 0$, so for all $m\in\bN$, 
$$
\nm {(UT)^m}\ge |\lam_\bc(T)|^m\cdot 2^{-L(1+\rho(t))}= |\lam_\bc(T)|^m\cdot 2^{-m(1+l)(1+\rho(t))}.
$$
  Accordingly $UT\in \abws$ is not a quasinilpotent operator, and $T\notin\rad \abws$. \endproof

\section{$\cA^{**}$ is semisimple.}

   Let $\theta_0:\cT^{***}\to\cT^*=B(\cH)$ be the natural projection, which is an algebra homomorphism, and let $\theta=\theta_0|_{\cA^{**}}$ be the restriction, which is a representation of $\cA^{**}$. If $\tau\in\cA^{**}$ is a weak-* limit of operators $T_\alp$ in $\cA$, then for each 
$\eta,\zeta\in\cH$, we have $\sprod{\theta(\tau)\eta}\zeta=\lim_\alp\sprod{T_\alp\eta}{\zeta}$, so $\theta(\tau)$ is the $\sig(B(\cH),\cT)$-limit of the operators $T_\alp$, and the image $\theta(\cA^{**})$ is the weak-* closure $\abws$ of $\cA$ in $B(\cH)$. Conversely, the image of the unit ball of $\cA^{**}$, being the weak-* continuous image of a weak-* compact set, is weak-* compact, and therefore contains the weak-* closure $\bar B^{w*}$ of the unit ball of $\cA$.  It is a consequence of the Hahn-Banach theorem that $\abws$ is equal to the union $\cupp n1\infty n\cdot    \bar B^{w*}$, so 
we have $\theta(\cA^{**})=\abws$, which by \refthm{t4}\ is semisimple. To deduce that $\cA^{**}$ is semisimple, we need only prove that $\theta$ is a faithful (injective) representation. 

\bthm\label{t5} The representation $\theta:\cA^{**}\to B(H)$ is faithful.
\ethm

\proof Let $\tau\in  \cA^{**}$ with $\nm\tau=1$. We claim that $\theta(\tau)\ne 0$. To establish this, we first prove the following Lemma: 

\blem\label{5-1} If $\tau\in  \cA^{**}$ with $\nm\tau=1$, then for every $\veps>0$, there is an $n\in\bN$ and a $\phi\in B(H)^*$ with $\nm\phi=1$, such that the 
compression $\phi_n=\bar P_n\cdot\phi\cdot\bar P_n$ satisfies $|\sprod{\tau}{\phi_n}|>1-\veps$. 
\elem 
\proof If $a,b\in B(H)$ and $Q$ is an orthogonal projection, then simple calculations yield the  inequalities
$\nm{aQ+b(1-Q)}\le\sqrt{\nm{aQ}^2+\nm{b(1-Q)}^2}$ and $\nm{Qa+(1-Q)b}\le\sqrt{\nm{Qa}^2+\nm{(1-Q)b}^2}$. When these are dualized, 
the directions of the inequalities are reversed: if $\phi,\psi\in B(H)^*$ then
$$
\nm{\phi\cdot Q+\psi\cdot(1-Q)}\ge\sqrt{\nm{\phi\cdot Q}^2+\nm{\psi\cdot(1-Q)}^2},
$$
and
\beq\label{27}
\nm{Q\cdot\phi+(1-Q)\cdot\psi}\ge\sqrt{\nm{Q\cdot\phi}^2+\nm{(1-Q)\cdot\psi}^2}. 
\eeq
For every $\eta>0$ there is a $\phi\in \cA^*$ such that $\nm\phi=1$ and $\sprod{\tau}{\phi}>1-\eta$. There is also a witness 
$T\in\cA$ such that $\nm T=1$ and $\sprod{\phi}T>1-\eta$. By \reflem{l3}\ part (c), there is an $n\in\bN$ such that $\nm{(1-\bar P_n)T}<\eta$. 
Hence, $|\sprod {\phi- \phi_n}T|\le \nm{(1-\bar P_n)T}+\nm{\bar P_nT(1-\bar P_n)}=\nm{(1-\bar P_n)T}$ (because $\ker\bar P_n$ is an invariant subspace for $\cA$) $< \eta$ also, and so $\nm{\phi_n}\ge |\sprod {\phi_n}T|> 1-2\eta$. By \refeq{27}\ we therefore have   $\nm{(1-\bar P_n)\cdot\phi },\nm{\phi\cdot (1-\bar P_n)}<\sqrt{1-(1-2\eta)^2}<2\sqrt{\eta}$, hence $\nm{\phi-\phi_n}<4\sqrt\eta$. Since $\sprod{\tau}{\phi}>1-\eta$, we have 
$|\sprod \tau{\phi_n}|> 1-\eta-\nm{\phi-\phi_n}\ge 1-\eta-4\sqrt{\eta}$. Appropriate choice of $\eta>0$ yields $ |\sprod \tau{\phi_n}|>1-\veps$ as required.
\endproof

We now prove \refthm{t5}.
Let $\tau\in\cA^{**}$ with $\nm\tau=1$, and assume towards a contradiction that $\theta(\tau)=0$.
Write $\gam_n=\sup\{|\sprod{\bar P_n\cdot\phi\cdot\bar P_n}\tau|:\phi\in B(H)^*,\nm\phi=1\}$.
The sequence $\gam_n$ is non-decreasing, and by \reflem{5-1}\ we have $\gam_n\to 1$. Pick then an $N\in\bN$ such that 
$\gam_N>0$, and let $n\le N$ be the least natural number such that  $\gam_n=\gam_N$. 
For every $\veps>0$ we can find $\phi\in B(H)^*$, $\nm\phi=1$ such that $\sprod{\bar P_n\cdot\phi\cdot\bar P_n}\tau\ge\gam_N-\veps$. 

Given such an $\veps>0$ and  $\phi$, we write $\phi_1$ for a weak-* accumulation point of the functionals $\bar\pi_k\cdot \phi$; but actually, we claim that  $\phi_1$ is the norm convergent limit of $\bar\pi_k\cdot \phi$. For the norms 
$\nm{\bar\pi_k\cdot \phi}$ are a nondecreasing sequence tending to a limit $l$; equation \refeq{27}\ tells us that for $m>k$ we have $\nm{\bar\pi_k\cdot \phi}^2+\nm{(\bar\pi_m-\bar\pi_k)\cdot \phi}^2\ge\nm{\bar\pi_m\cdot \phi}^2$, so 
$\nm{(\bar\pi_m-\bar\pi_k)\cdot \phi}\to 0$ as $k,m\to\infty$; so the sequence $(\bar\pi_k\cdot \phi)_{k\in\bN}$ satisfies the Cauchy criterion and is norm convergent. Each projection ${\bar\pi_k}$ is of finite rank, so $\bar\pi_k\cdot \phi$ belongs to the trace-class operators $\cT$. Therefore, $\phi_1\in \cT$. But the difference $\phi-\phi_1=\lim_k(1-\bar\pi_k)\cdot\phi$ will annihilate any compact operator.

We therefore claim that $n>1$. For by \refthm{t2}, whenever $T\in\cA$ the operator $P_1T=P_1T\up 1$ is a compact operator, so $\sprod {P_1TP_1}\phi=\sprod {P_1TP_1}{\phi_1}$. We may write $\tau$ as a weak-* convergent limit $\tau=\lim_{w*}T_\alp$ for $T_\alp\in\cA$ with $\nm{T_\alp}=1$. Then 
$\gam_N-\veps\le\sprod {P_1\cdot\phi\cdot P_1}\tau=\lim_\alp \sprod {T_\alp}{P_1\cdot\phi\cdot P_1}=\lim_\alp
\sprod {P_1T_\alp P_1}\phi=\lim_\alp\sprod {P_1T_\alp P_1}{\phi_1}=\lim_\alp\sprod {T_\alp}{P_1\cdot\phi_1\cdot P_1}
 =\sprod {\tau}{P_1\cdot\phi_1\cdot P_1}$. For small $\veps$ this implies $\sprod {\tau}{P_1\cdot\phi_1\cdot P_1}\ne 0$. But
$P_1\cdot\phi_1\cdot P_1\in \cT=B(H)_*$, so $\theta(\tau)$ is not the zero operator in $B(H)$, a contradiction. Therefore we have $n>1$. 

Given $n>1$, we again pick $\veps>0$ and find $\phi\in B(H)^*$, $\nm\phi=1$ such that 
\beq\label{28}
\sprod{\bar P_n\cdot\phi\cdot\bar P_n}\tau\ge\gam_N-\veps>0.
\eeq
 The norm limit $\phi_1=\lim_k\bar\pi_k\cdot \phi$ is again in 
$\cT$. However, the difference $\phi-\phi_1$ will not necessarily annihilate $\bar P_n T\bar P_n$ for $T\in\cA$, because
though $\phi-\phi_1$ annihilates $K(H)$, 
the operator $\bar P_n T\bar P_n$ need not be compact. Rather, for $T\in\cA$ we have $\bar P_n T\bar P_n=\bar P_n\bar T\up n\bar P_n$, where $\bar T\up n=\bar\Del_n(T)$ as in \refdfn{d3b}; and $\bar T\up n=\bar T\up {n-1}+T\up n$, where the operator 
$\bar P_nT\up n$ is compact by \refthm{t2}. So $\sprod  {\bar P_nT\up n\bar P_n}{\phi-\phi_1}=0$ for all $T\in\cA$. Writing 
$\tau=\lim_\alp T_\alp$ for a suitable net $(T_\alp)$ in $\cA$, we have  
$\sprod  {T_\alp\up n}{\bar P_n(\phi-\phi_1)\bar P_n}=0$ for all $\alp$. Writing $\beta=\lim_\alp\sprod  {\bar T_\alp\up {n-1}}{\bar P_n(\phi-\phi_1)\bar P_n}$, we will have 
$0=\sprod{\bar P_n\phi_1\bar P_n}\tau$ (because $\phi_1\in\cT$ and $\theta(\tau)=0$ by hypothesis) $=\sprod{\bar P_n\phi\bar P_n}\tau-\lim_\alp\sprod{\bar P_n(\phi-\phi_1)\bar P_n}{T_\alp}
=$
$\sprod{\bar P_n\phi\bar P_n}\tau-\beta$. By equation \refeq{28}, we have $|\bet|\ge \gam_N-\veps$. 

For each $T\in\cA$ and $n>1$, the norms of $\bar T\up {n-1}$ and $\bar P_{n-1}\bar T\up {n-1}\bar P_{n-1}=\bar P_{n-1}T\bar P_{n-1}$ are the same by \refeq{15}. Thus there is a unique map $\eta:\bar P_{n-1}\cdot\cA\cdot \bar P_{n-1}\to\bar\cA\up{n-1}$ which is a right inverse to the compression $p:\cA\to \bar P_{n-1}\cdot\cA\cdot \bar P_{n-1}$ with $p(T)=\bar P_{n-1}T\bar P_{n-1}$ ($T\in\cA$); and $\nm\eta=1$. We will have $\eta\cdot p=\bar\Del_{n-1}$. Let us write $\psi=(\bar P_n(\phi-\phi_1)\bar P_n)\circ \eta$. Then 
$\psi\in (\bar P_{n-1}\cdot\cA\cdot \bar P_{n-1})^*$ with $\nm\psi\le 1$. 

By the Hahn-Banach theorem, we can extend $\psi$ to $\bar P_{n-1}\cdot B(H)\cdot \bar P_{n-1}$ with the same norm; and then extend to all of $B(H)$ so that $\psi=\psi\circ p$ (where we abuse notation slightly by writing $p$ for the compression $B(H)\to \bar P_{n-1}\cdot B(H)\cdot \bar P_{n-1}$ also).

Then for $T\in\cA$ we have $\psi(T)=\psi\circ \eta p(T)=\psi(\bar\Del_{n-1}(T))$; so 
$$
|\sprod\psi\tau|=\lim_\alp |\sprod \psi{ \bar \Del_{n-1}T_\alp}|
=\lim_\alp|\sprod{\bar P_n(\phi-\phi_1)\bar P_n}{\bar T_\alp\up{n-1}}|=|\bet|\ge\gam_N-\veps.
$$
Since $\psi=\psi\circ p=\bar P_{n-1}\cdot\psi\cdot\bar P_{n-1}$, we find that $\gam_{n-1}=\sup\{|\sprod{\bar P_{n-1}\cdot\phi\cdot\bar P_{n-1}}\tau|:\phi\in B(H)^*,\nm\phi=1\}$ is at least $\gam_N-\veps$. But $\veps>0$ is arbitrary, so $\gam_{n-1}=\gam_N$, and $n$ was not the minimal integer with 
$\gam_n=\gam_N$ contrary to hypothesis. This contradiction proves the Theorem.\endproof

\section{References to Gulick's paper.}

Having established that the result $\rad A^{**}\cap A=\rad A$ of Gulick is wrong, let us look at papers which have referenced Gulick \cite{G} since 1966, and try to establish that no further damage has been done.

The lengthy paper of Dales and Lau \cite{DL} refers to Gulick's paper \cite{G}, but does not use the false theorem 4.6; private communication with my colleague Garth Dales reveals a history of previous suspicion about that result, but no actual counterexamples as presented here.  The paper of Daws, Haydon, Schlumprecht and White \cite{DHSW} refers to (the proof of) Theorem 3.3 of \cite{G}, which we believe to be completely correct. Likewise the paper of Bouziad and Filali \cite{BF} quotes the proof, given by Gulick in \cite{G} (Lemma 5.2), that the radical of $L^\infty(G)^*$ is nonseparable for any non-discrete locally compact group $G$. This proof also is perfectly valid. The earlier paper of Granirer \cite{Gr}\ makes reference to that same, correct, Lemma. Tomiuk \cite{T} likewise refers to Gulick's untainted Theorem 5.5. 
In \cite{U}, A. \"{U}lger solves one of the problems posed by Gulick in \cite{G}. Finally Tomiuk and Wong \cite{TW}\ make a passing reference to \cite{G} in their paper on Arens products.  

We have not found a case in which another author has used the false Theorem 4.6 from Gulick's paper, or anything tainted by it. This chimes with our reckoning that more than one author apart from ourselves has suspected that that Theorem is false. So, the general literature on Banach algebras is not seriously harmed; but it was nonetheless high time that these counterexamples were made known so that such errors will not occur in the future.

\end{document}